\documentclass{scrartcl}

\usepackage{amsmath}
\usepackage{amsfonts}
\usepackage{amssymb}

\usepackage{tikz}

\usepackage{hyperref}

\title{Foundations of brick diagrams}
\author{Jules Hedges \and Jelle Herold}

\begin{document}

\maketitle

\begin{abstract}
	We discuss the foundations of 2-dimensional graphical languages, with a view towards their computer implementation in a `compiler' for monoidal categories.
	In particular, we discuss the close relationship between string diagrams, pasting diagrams, linear logic proof trees and k-d trees, the last being a data structure from computational geometry and computer graphics that recursively partitions a space.
	We introduce a minor variant of cubical pasting diagrams, which we call `brick diagrams', which are used in the Statebox visual programming language.
	This paper is intended as a discussion and literature review, and contains little mathematics.
\end{abstract}

\section{Introduction}

There has been much talk in the fledgeling Applied Category Theory community of the importance of a general-purpose \emph{compiler} for string diagrams, distinct from the several \emph{proof assistants} in development.
The purpose of such a tool is first to assist the user in building and visualising a string diagram, and second to automate the process of computing the denotation in a semantic category.
For example, a string diagram in the prop of matrices denotes a matrix, and the software should explicitly compute this matrix, or produce code in (for example) NumPy to do so.
In this paper we propose a foundation for such a tool, based on a four-way connection between \emph{string diagrams}, \emph{cubical pasting diagrams}, \emph{proof trees} and \emph{k-d trees}.

There are, broadly, two ways of presenting free monoidal categories.
The first is categories of string diagrams modulo isotopy, and the second is as categories of linear logic proof trees modulo rewrites.
The former is far more useful for calculating with, but the second is finitistic and much easier to implement in a computer.
We propose a bridge between the two via \emph{k-d trees}, a data structure from computer graphics that can be used to represent string diagrams in a way that coincides with proof trees.

The importance of string diagrams can be expressed formally by the fact (better known as the Joyal-Street coherence theorem \cite{joyal91}) that they are a presentation of \emph{free monoidal categories}.
More precisely, the free monoidal category on a collection of generating objects and a collection of generating morphisms, is equivalent to a category whose morphisms are isotopy classes of string diagrams with strings labelled by generating objects and nodes labelled by generating morphisms.
The universal property of freeness states that if each generating object and generating morphism is associated to an object and morphism in any target `semantic' monoidal category, then any diagram induces a canonical morphism, its `semantics', which is isotopy-invariant.
This property is so intuitively obvious that it is still used routinely by many people who do not fully understand the formalisation.

(Side note: In general, freeness is a property of \emph{syntax}, and string diagrams are a 2-dimensional syntax.
Decategorifying, the free \emph{monoid} on a set of generating symbols is the monoid of \emph{strings} of symbols (also known as lists), which are used as syntax in ordinary (1-dimensional) algebra.)

By instantiating with an example, this becomes much less abstract.
Suppose that our target category is $\mathbf{Mat}_\mathbb R$, the category whose objects are natural numbers and whose morphisms $m \to n$ are $m \times n$ matrices, with matrix multiplication as composition.
Suppose we take a single generating object $1$, and take each $m \times n$ matrix as a generating morphism of type $1^{\otimes m} \to 1^{\otimes n}$.
The Joyal-Street coherence theorem then states that any string diagram in which every node is labelled by a matrix (whose dimensions match the node's input and output degree) determines a matrix, which is isotopy-invariant.

There is moreover a simple procedure to calculate this matrix.
We recursively cut the diagram horizontally and vertically, and stop when each region is indivisible, containing either a single node, or a single string, or a single string crossing.
Then we associate each indivisible region to the associated semantics (respectively either the labelling matrix, or the $1 \times 1$ identity matrix $\begin{pmatrix} 1 \end{pmatrix}$, or the swap matrix $\begin{pmatrix} 0 & 1 \\ 1 & 0 \end{pmatrix}$), and combine them using matrix multiplication and Kronecker product according to the horizontal and vertical cuts we made.
The Joyal-Street theorem states that the resulting matrix is independent of our choices of cuts.
Moreover, although we illustrated the procedure with matrices, the procedure is the same in any monoidal category.

By a \emph{compiler} for monoidal categories we mean a piece of software that will carry out this process automatically.

\subsection{Computer implementation of string diagrams}

The need for computer assistance has been widely discussed in the fledgeling Applied Category Theory community, for example at the first ACT conference in Leiden in 2018, the second in Oxford in 2019, and on the ACT Google Groups mailing list.

There are several existing tools for working with string diagrams that are classified as \emph{proof assistants}.
Quantomatic and pyZX are tools for simplifying diagrams in the language of ZX calculus, which internally use a graphical representation (i.e. based on graphs) and are essentially specialised graph rewrite engines.
By using a graphical representation, these tools are essentially specialised to compact closed categories, if not to ZX calculus specifically.

Homotopy.io and its predecessor Globular are tools for constructing morphisms in higher categories using their associated diagrams.
Homotopy.io is based on an entirely new definition of higher categories that was partly designed with computer implementation in mind \cite{dorn_associative_n_categories}.
However we argue that while it is clearly an effective foundation for a proof assistant it does not seem to be well-suited to compilation, and that gains (in efficiency and mathematical difficulty) can still be made by focussing on the 2-dimensional case.

What is conspicuously absent is any tool that can be classified as a \emph{compiler}, for working with the semantics rather than the syntax of string diagrams.
Consider, for example, the prop $\mathbf{Mat}_\mathbb R$ whose objects are natural numbers and whose morphisms $m \to n$ are $m \times n$ real matrices.
With such a tool, the user would create a diagram using a graphical user interface in which every atomic block is labelled by a matrix of the appropriate dimensions.
The user will then push a button and, if the diagram is well-formed, the software will compute the matrix represented by the diagram (or generate code that will do so, for example in numPy).

In general, (applied) category theory can make connections between disparate application areas by finding commonalities between their graphical languages.
A common front end user interface that is able to support only a small number of topological flavours of string diagrams (for example monoidal, symmetric monoidal, traced monoidal, compact closed) together with customisable generating morphisms will be applicable to many different domain-specific backends.
In each case, the backend must supply a `meaning' to the morphisms as a representation in the computer, that supports the operations needed to interpret the relevant class of string diagrams (chiefly, sequential and parallel composition operators).
For example, in the case of matrix algebra, each morphism represents a matrix, and sequential and parallel composition correspond semantically to matrix multiplication and Kronecker product.

There is a relatively obvious way to organise software like this, which many people have thought of independently.
The front end will generate `code' in an \emph{intermediate language}, which is simply a representation of a term in the logical language of monoidal categories, implemented in a standard meta-language such as XML or JSON.
Then each application domain can be supported by a separate backend compiler that translates this intermediate language into code or a file format that can be read by a domain-specific application, for example Simulink.
Each of these components is decoupled and communicates only through the intermediate language, which is sociologically important because each component can be written in a programming language familiar to its domain community.

A diagram can generally be `compiled' in many different ways, corresponding to different ways of `carving up' the diagram using horizontal and vertical cuts, and also by doing topological moves on the diagram.
The axioms of a monoidal category (for example) state precisely that all of these different ways will produce the same end result.
However, they may be vastly different from a computational perspective.
Taking the example of matrices again, Kronecker product is exponentially more expensive than matrix multiplication, and so a good heuristic will be to manipulate the diagram in order to minimise the maximum number of parallel compositions, aka. $\otimes$-width.

In the end, any general-purpose software will need to be scriptable to allow domain-specific heuristics.
As a starting point, however, any compilation method will be suitable for small diagrams in every application domain.

\section{Graphs}

Probably the most obvious way to represent a string diagram is as a graph, which we will call a \emph{graphical} representation.
Of course, the representation of graphs in a computer is extremely well-studied with efficient implementations available in every serious programming language.
Probably the most common representation is \emph{adjacency lists}, where the graph is represented by a list of vertices, where each vertex is associated to a list of adjacent vertices.
(Note however that representation of graphs in purely functional programming languages is notoriously tricky, with the Haskell's Algebraic Graph Library \cite{mokhov_algebraic_graphs} being a recent development.)

Graphical representations are used by the Quantomatic \cite{kissinger_zamdzhiev_quantomatic} and PyZX \cite{kissinger_wetering_pyzx} proof assistants for ZX calculus, as well as Evan Patterson's recent Catlab library\footnote{\url{https://github.com/epatters/Catlab.jl}} for the scientific programming language Julia.

Graphs are certainly reasonable for use-cases that are specialised to compact closed categories (as the above examples are), since any node may be connected to any other \cite{kelly_laplaza_coherence_compact_closed_categories}.
However, it is also easy to generalise from compact closed to symmetric monoidal categories, by restricting to directed acyclic graphs.

Certainly there are many examples of compact closed categories of importance in applied category theory.
Examples include the ZX calculus used in categorical quantum mechanics \cite{coecke_kissinger_picturing_quantum_processes}, the category of finite dimensional vector spaces used in distributional semantics \cite{coecke_sadrzadeh_clarke_discocat}, and categories of decorated cospans \cite{fong_decorated_cospans} and corelations \cite{fong_decorated_corelations}.

However, while compact closure is a `gold standard' for categorical structure, there are also important monoidal categories that are not compact closed, for example those used in dataflow \cite{delpeuch_complete_language_faceted_dataflow_programs}, classical computing \cite{pavlovic13} (including the Statebox visual programming language \cite{statebox_monograph}), and compositional game theory \cite{hedges_etal_compositional_game_theory}.
We note however that all of these examples are symmetric monoidal categories.
(Non-symmetric monoidal categories find application for example in linguistics, topology and parts of theoretical physics, but these are certainly the exception.)

The formal definition of string diagrams \cite{joyal91} is actually via \emph{topological graphs}, which are Hausdorff spaces that are 1-dimensional manifolds plus singular points representing the nodes (which is equivalently a 1-dimensional simplicial complex \cite[section 1.A]{hatcher_algebraic_topology}), equipped with an embedding in the plane.
(Notice that an embedding is by definition injective, which forces planarity.)

An important point to keep in mind is that when drawing string diagrams in practice we nearly always abuse notation, and draw the singular points as regions (of various shapes, which can themselves carry information, such as when a morphism is a transpose or adjoint \cite[section 4.3.2]{coecke_kissinger_picturing_quantum_processes}).
For example, instead of the graph
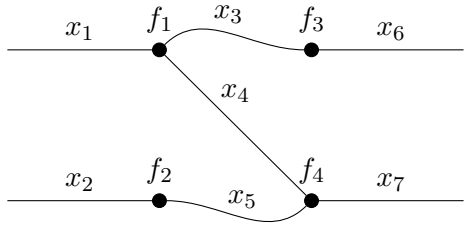
\begin{center} \begin{tikzpicture}
	\node (f1) [circle, scale=.5, fill=black, label={$f_1$}, draw] at (2, 2) {}; \node (f2) [circle, scale=.5, fill=black, label={$f_2$}, draw] at (2, 0) {};
	\node (f3) [circle, scale=.5, fill=black, label={$f_3$}, draw] at (4, 2) {}; \node (f4) [circle, scale=.5, fill=black, label={$f_4$}, draw] at (4, 0) {};
	\draw [-] (0, 2) to node [above] {$x_1$} (f1); \draw [-] (0, 0) to node [above] {$x_2$} (f2);
	\draw [-] (f1) to [out=45, in=180] node [above] {$x_3$} (f3);
	\draw [-] (f1) to [out=-45, in=135] node [above=5pt] {$x_4$} (f4);
	\draw [-] (f2) to [out=0, in=-135] node [above] {$x_5$} (f4);
	\draw [-] (f3) to node [above] {$x_6$} (6, 2); \draw [-] (f4) to node [above] {$x_7$} (6, 0);
\end{tikzpicture} \end{center}
we would typically draw
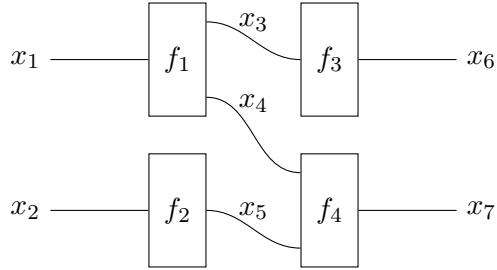
\begin{center} \begin{tikzpicture}
	\node (x1) at (0, 2) {$x_1$}; \node (x2) at (0, 0) {$x_2$}; \node (x6) at (6, 2) {$x_6$}; \node (x7) at (6, 0) {$x_7$};
	\node (f1) [rectangle, minimum height=1.5cm, minimum width=.75cm, draw] at (2, 2) {$f_1$};
	\node (f2) [rectangle, minimum height=1.5cm, minimum width=.75cm, draw] at (2, 0) {$f_2$};
	\node (f3) [rectangle, minimum height=1.5cm, minimum width=.75cm, draw] at (4, 2) {$f_3$};
	\node (f4) [rectangle, minimum height=1.5cm, minimum width=.75cm, draw] at (4, 0) {$f_4$};
	\node (dummy1) at (0, 2.5) {}; \node (dummy2) at (0, 1.5) {}; \node (dummy3) at (0, .5) {}; \node (dummy4) at (0, -.5) {};
	\draw [-] (x1) to (f1); \draw [-] (x2) to (f2);
	\draw [-] (f1.east |- dummy1) to [out=0, in=180] node [above] {$x_3$} (f3);
	\draw [-] (f1.east |- dummy2) to [out=0, in=180] node [above=5pt] {$x_4$} (f4.west |- dummy3);
	\draw [-] (f2) to [out=0, in=180] node [above] {$x_5$} (f4.west |- dummy4);
	\draw [-] (f3) to (x6); \draw [-] (f4) to (x7);
\end{tikzpicture} \end{center}
(We will continue to use this as a running example in the next few sections.)

This can be thought of as infinitesimal information.
To our knowledge, this has not been formalised using differential topology.

\section{Planar graphs and cell complexes}

Given that most of the categories we care about are symmetric, why should we care about planar graphs?
The answer is that the symmetry of a monoidal category can be computationally nontrivial, and we may want finer control over it.
So, even if our target category is symmetric, we may want to represent its diagrams as planar graphs together with a special type of node representing a string crossing.
(For example, in the prop of matrices a string crossing denotes the swap matrix $\begin{pmatrix} 0 & 1 \\ 1 & 0 \end{pmatrix}$.)
In addition to this, planar string diagrams are slightly more general; on the other hand, it would also be possible to begin with a graphical representation and later add the option of enforcing planarity.

There has been some work, including some very recent, on data structures for efficient editing of planar graphs \cite{tamassia_online_planar_graph_embedding,holm_rotenberg_dynamic_planar_embeddings,karczmarz_data_structures_dynamic_graphs}.

A more obvious representation of planar graphs is using \emph{polygonal cell complexes}, i.e. the 2-dimensional case of \emph{polytopic cell complexes}.
(This implicitly uses F\'ary's theorem, which says that every planar graph has a planar embedding in which all edges are straight lines.)
This means that the focus is on the regions of the graph, where each region has associated to it a pair of lists of the nodes around its edges, in order.

While detecting planarity of an abstract polygonal complex may be no easier than detecting planarity of a graph, the advantage of cell complexes is that, given an existing polygonal cell complex that is planar, it is essentially trivial to decide when a new edge can be added in a way that preserves planarity: namely, when it is between two vertices bordering the same cell.

\section{Pasting diagrams}

The Poincar\'e dual of a string diagram is a \emph{pasting diagram}.
There are different flavours of pasting diagram, corresponding to different flavours of 2-category.

If we view our strict monoidal category as a 2-category with 1 object, we can also represent our example morphism by a \emph{globular pasting diagram}:
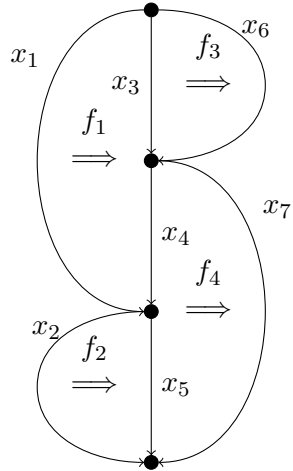
\begin{center} \begin{tikzpicture}
	\node (n1) [circle, scale=.5, fill=black, draw] at (0, 6) {}; \node (n2) [circle, scale=.5, fill=black, draw] at (0, 4) {};
	\node (n3) [circle, scale=.5, fill=black, draw] at (0, 2) {}; \node (n4) [circle, scale=.5, fill=black, draw] at (0, 0) {};
	\draw [->] (n1) to [out=180, in=90] node [left=5pt] {$x_1$} (-1.5, 4) to [out=-90, in=180] (n3);
	\draw [->] (n3) to [out=180, in=90] node [left] {$x_2$} (-1.5, 1) to [out=-90, in=180] (n4);
	\draw [->] (n1) to node [left] {$x_3$} (n2); \draw [->] (n2) to node [right] {$x_4$} (n3); \draw [->] (n3) to node [right] {$x_5$} (n4);
	\draw [->] (n1) to [out=0, in=90] node [right] {$x_6$} (1.5, 5) to [out=-90, in=0] (n2);
	\draw [->] (n2) to [out=0, in=90] node [right=5pt] {$x_7$} (1.5, 2) to [out=-90, in=0] (n4);
	\node [label={$f_1$}] at (-.75, 4) {$\implies$}; \node [label={$f_2$}] at (-.75, 1) {$\implies$};
	\node [label={$f_3$}] at (.75, 5) {$\implies$}; \node [label={$f_4$}] at (.75, 2) {$\implies$};
\end{tikzpicture} \end{center}

Globular pasting diagrams are unwieldy to work with for any example much larger than this one, and in the authors' opinions they are also less intuitive in practice, in spite of the fact that it seems more intuitive that 2-cells go between 1-cells just as morphisms go between objects.
Indeed, the introduction of string diagrams over the older pasting diagrams is intimately connected with the emergence of modern applied category theory.

However, a strict monoidal category can also be viewed not only as a 2-category with one object, but also as a double category with one 0-cell and one horizontal 1-cell.
These \emph{cubical pasting diagrams} overcome the practical shortcomings of globular pasting diagrams with respect to actually drawing them on a page:
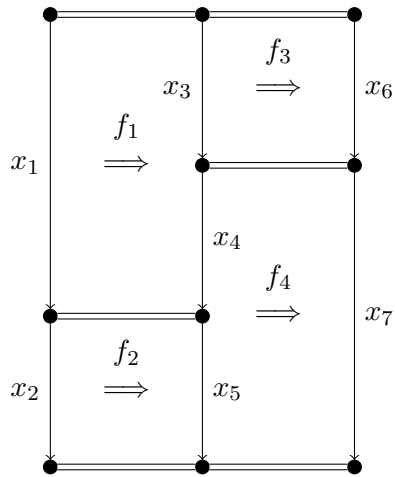
\begin{center} \begin{tikzpicture}
	\node (n1) [circle, scale=.5, fill=black, draw] at (0, 6) {}; \node (n2) [circle, scale=.5, fill=black, draw] at (0, 2) {}; \node (n3) [circle, scale=.5, fill=black, draw] at (0, 0) {};
	\node (n4) [circle, scale=.5, fill=black, draw] at (2, 6) {}; \node (n5) [circle, scale=.5, fill=black, draw] at (2, 4) {};
	\node (n6) [circle, scale=.5, fill=black, draw] at (2, 2) {}; \node (n7) [circle, scale=.5, fill=black, draw] at (2, 0) {};
	\node (n8) [circle, scale=.5, fill=black, draw] at (4, 6) {}; \node (n9) [circle, scale=.5, fill=black, draw] at (4, 4) {}; \node (n10) [circle, scale=.5, fill=black, draw] at (4, 0) {};
	\draw [double equal sign distance] (n1) to (n4); \draw [double equal sign distance] (n4) to (n8);
	\draw [double equal sign distance] (n2) to (n6); \draw [double equal sign distance] (n5) to (n9);
	\draw [double equal sign distance] (n3) to (n7); \draw [double equal sign distance] (n7) to (n10);
	\draw [->] (n1) to node [left] {$x_1$} (n2); \draw [->] (n2) to node [left] {$x_2$} (n3);
	\draw [->] (n4) to node [left] {$x_3$} (n5); \draw [->] (n5) to node [right] {$x_4$} (n6); \draw [->] (n6) to node [right] {$x_5$} (n7);
	\draw [->] (n8) to node [right] {$x_6$} (n9); \draw [->] (n9) to node [right] {$x_7$} (n10);
	\node [label={$f_1$}] at (1, 4) {$\implies$}; \node [label={$f_2$}] at (1, 1) {$\implies$};
	\node [label={$f_3$}] at (3, 5) {$\implies$}; \node [label={$f_4$}] at (3, 2) {$\implies$};
\end{tikzpicture} \end{center}

As a further piece of evidence in favour of cubical diagrams being more intuitive, we note that Joyal and Street's original coherence proof for monoidal categories \cite{joyal91} proceeds by division of the plane into rectangular regions, even though globular regions would have sufficed.

Pasting diagrams of various shapes have been formalised in terms of \emph{pasting schemes} \cite{johnson_pasting_diagrams,power_pasting_theorem} and \emph{parity complexes} \cite{street_parity_complexes}.

\section{Tilings and brick diagrams}

We make some minor notational changes to the previous cubical pasting diagram, to make it visually cleaner:
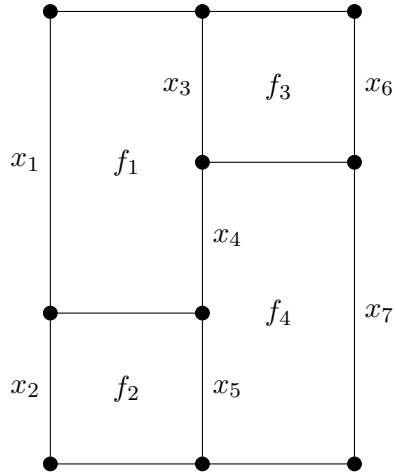
\begin{center} \begin{tikzpicture}
	\node (n1) [circle, scale=.5, fill=black, draw] at (0, 6) {}; \node (n2) [circle, scale=.5, fill=black, draw] at (0, 2) {}; \node (n3) [circle, scale=.5, fill=black, draw] at (0, 0) {};
	\node (n4) [circle, scale=.5, fill=black, draw] at (2, 6) {}; \node (n5) [circle, scale=.5, fill=black, draw] at (2, 4) {};
	\node (n6) [circle, scale=.5, fill=black, draw] at (2, 2) {}; \node (n7) [circle, scale=.5, fill=black, draw] at (2, 0) {};
	\node (n8) [circle, scale=.5, fill=black, draw] at (4, 6) {}; \node (n9) [circle, scale=.5, fill=black, draw] at (4, 4) {}; \node (n10) [circle, scale=.5, fill=black, draw] at (4, 0) {};
	\draw [-] (n1) to (n4); \draw [-] (n4) to (n8);
	\draw [-] (n2) to (n6); \draw [-] (n5) to (n9);
	\draw [-] (n3) to (n7); \draw [-] (n7) to (n10);
	\draw [-] (n1) to node [left] {$x_1$} (n2); \draw [-] (n2) to node [left] {$x_2$} (n3);
	\draw [-] (n4) to node [left] {$x_3$} (n5); \draw [-] (n5) to node [right] {$x_4$} (n6); \draw [-] (n6) to node [right] {$x_5$} (n7);
	\draw [-] (n8) to node [right] {$x_6$} (n9); \draw [-] (n9) to node [right] {$x_7$} (n10);
	\node at (1, 4) {$f_1$}; \node at (1, 1) {$f_2$};
	\node at (3, 5) {$f_3$}; \node at (3, 2) {$f_4$};
\end{tikzpicture} \end{center}

Geometrically, this is a \emph{rectangular tiling} of a rectangle.
These rectangular tilings have been specifically studied, with \cite{dawson_pare_characterising_tileorders} giving combinatorial representations in terms of \emph{double orders}, as well as a small generating set under a `welding' operation.

An issue with this representation is that we need to treat the monoidal unit as a special case, as a labelled edge that can always be introduced or cancelled (or give up the idea of using rectangular tilings, which will lead to new difficulties).
We propose the following notation, in which we take an additional Poincar\'e dual of only the vertical boundaries:
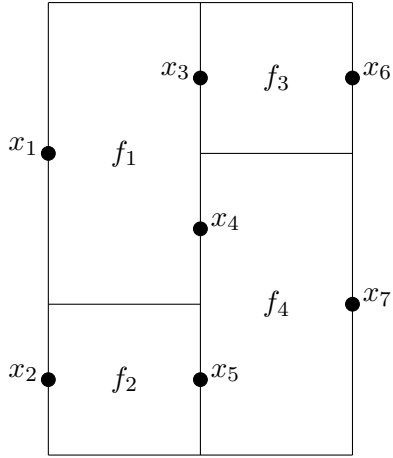
\begin{center} \begin{tikzpicture}
	\draw [-] (0, 6) to (4, 6); \draw [-] (0, 0) to (4, 0); \draw [-] (0, 6) to (0, 0); \draw [-] (4, 6) to (4, 0);
	\draw [-] (2, 6) to (2, 0); \draw [-] (0, 2) to (2, 2); \draw [-] (2, 4) to (4, 4);
	\node [circle, scale=.5, fill=black, label={[left] $x_1$}, draw] at (0, 4) {}; \node [circle, scale=.5, fill=black, label={[left] $x_2$}, draw] at (0, 1) {};
	\node [circle, scale=.5, fill=black, label={[left] $x_3$}, draw] at (2, 5) {};
	\node [circle, scale=.5, fill=black, label={[right] $x_4$}, draw] at (2, 3) {};
	\node [circle, scale=.5, fill=black, label={[right] $x_5$}, draw] at (2, 1) {};
	\node [circle, scale=.5, fill=black, label={[right] $x_6$}, draw] at (4, 5) {}; \node [circle, scale=.5, fill=black, label={[right] $x_7$}, draw] at (4, 2) {};
	\node at (1, 4) {$f_1$}; \node at (1, 1) {$f_2$};
	\node at (3, 5) {$f_3$}; \node at (3, 2) {$f_4$};
\end{tikzpicture} \end{center}
We refer to this as a \emph{brick diagram}.

A tiling that can be built using sequential and parallel composition is called \emph{binarily composable}.
In \cite{dawson_forbidden_suborder_characterisation} it is shown that the only obstructions to a tiling being binarily composable are the \emph{pinwheel}
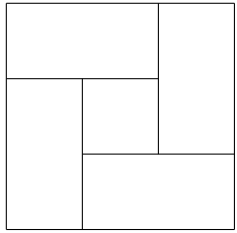
\begin{center} \begin{tikzpicture}
	\draw [-] (0, 3) to (3, 3); \draw [-] (0, 0) to (0, 3); \draw [-] (0, 0) to (3, 0); \draw [-] (3, 0) to (3, 3);
	\draw [-] (0, 2) to (2, 2); \draw [-] (2, 3) to (2, 1); \draw [-] (3, 1) to (1, 1); \draw [-] (1, 0) to (1, 2);
\end{tikzpicture} \end{center}
and its reflection.
If we restrict to binarily composable tilings then we can represent the tiling by its composition tree (which is usually not unique).
This is the subject of the next section.

Far more detail on tiling diagrams and pinwheels can be found in \cite{meyers_string_diagrams_double_categories}.
Note that for interpretation in a monoidal category (rather than a general double category) any pinwheel can be removed without changing the interpretation.
However, we still prefer our tilings to be pinwheel-free in order that we can represent them by their decompositions.

%\section{Tilings}
%
%A \emph{tile diagram} is a division of the square $[0, 1]^2$ into rectangles.
%We restrict to tile diagrams that do not contain a \emph{pinwheel}.
%In \cite{dawson_forbidden_suborder_characterisation} it is shown that this restriction is equivalent to being \emph{binarily composable}, which is to say, it can be described by an expression that builds it from atomic parts using sequential and parallel composition.
%(Such an expression is, of course, usually not unique.)
%
%We work in the context of a monoidal signature, i.e. a set $Ob(\Sigma)$ of \emph{object symbols}, a set $Mor(\Sigma)$ of \emph{morphism symbols}, and for every morphism symbol $f$ a pair of words $s(f), t(f)$ over $Ob (\Sigma)$, the \emph{source} and \emph{target} of $f$.
%
%A \emph{brick diagram} is a tile diagram together with:
%\begin{itemize}
%	\item To each minimal rectangle, an assigned morphism symbol, or the special symbol $id$
%	\item To each minimal vertical edge, an assigned word of object symbols
%	\item Such that for each minimal rectangle, concatenating together all the words along its entire left (right) edge yields its source (target)
%\end{itemize}
%
%Each minimal rectangle in a brick diagram is called a \emph{brick}, and should be thought of like a 2-dimensional Lego brick: its top and bottom edges have a pattern of studs and holes specified by its source and target, whereas its side edges are smooth.
%This means that brick diagrams are interpreted in monoidal categories rather than double categories; however, it should be straightforward to generalise to double categories.

\section{Proof trees and kd-trees}

Let $\Sigma$ be a monoidal signature.
We consider the noncommutative linear logic of $\otimes$ only, with elements of $Ob (\Sigma)$ as atomic propositions, and the following proof rules:
\[ (id) \frac{}{x \vdash x} (x \in Ob (\Sigma)) \qquad\qquad (f-ax) \frac{}{s (f) \vdash t(f)} (f \in Mor (\Sigma)) \]
\[ (\otimes-intro) \frac{a \vdash b \qquad a' \vdash b'}{aa' \vdash bb'} \qquad\qquad (cut) \frac{a \vdash b \qquad b \vdash c}{a \vdash c} \]

The free monoidal category on $\Sigma$ can be presented as the category whose objects are words over $Ob (\Sigma)$ and whose morphisms $a \to b$ are particular equivalence classes of proof trees with conclusion $a \vdash b$.
The equivalence relation that must be imposed is essentially generated by the axioms required of a monoidal category, known in proof theory as \emph{commuting conversions}.
(The authors have been unable to find a good reference for this fact; the general idea that proofs present morphisms in free categories is due to Lambek \cite{lambek_scott_introduction_higher_order_categorical_logic}.)

As a bridge between proof trees and diagrams, we consider \emph{binary space partition trees}, a common data structure in computational geometry \cite{berg_etal_computational_geometry}.
A binary space partition tree is a recursive partitioning of a space by dividing hyperplanes.
Specifically, we focus on \emph{k-d trees}, the special case in which the separating hyperplanes are all aligned to the coordinate axes.
We also focus on the 2-dimensional case.\footnote{Originally the term k-d tree was short for "k-dimensional tree", but standard terminology now is ``a k-d tree of dimension 2''.}

We consider k-d trees to store the nodes of a string diagram.
Thus we recursively partition the plane until each region contains only one node.
(It is also possible to consider a variant of k-d trees in which leaves are labelled by rectangular regions \cite{houthuys_box_sort}, in which we recover cubical pasting diagrams.)

K-d trees are not unique, and it seems obvious to conjecture that equivalences of k-d trees (which are well studied for the purpose of balancing) are the same as equivalences of proof trees.
(The same idea may also work in higher dimensions, but it is less obvious that nothing goes wrong, since higher dimensional diagrams are generally interpreted in semi-strict rather than strict $n$-categories.)

Now that we have introduced all the elements, we argue that the following are \emph{essentially the same}, modulo some unimportant details:
\begin{enumerate}
	\item Brick diagrams whose underling tiling is binarily composable
	\item K-d trees representing the nodes of the corresponding string diagram
	\item Proof trees for the above fragment of linear logic
\end{enumerate}
In particular, we argue that the appropriate equivalences that are imposed on each in order to present free monoidal categories are also essentially the same.

\section{Rendering string diagrams}

We do not expect many people (if any) to convert from working with string diagrams to tile diagrams day-to-day, and a software interface that forces tile diagrams on the user may put off many potential users.
Fortunately, we have the option to treat tiling diagrams only as a mental model and a data structure, and render a string diagram to the screen.
This can be done relatively easily by recursion on the k-d tree, performing a coordinate transformation at each recursive step.
This means, essentially, that the renderer is performing a Poincar\'e dualisation.

Although a language such as Javascript is probably better overall for implementing a graphical user interface, Haskell's \texttt{diagrams} package\footnote{\url{http://hackage.haskell.org/package/diagrams}} seems uniquely well-suited to this sort of operation.

At the base case of the recursion, rendering a region is done by drawing a node in its centre and then connecting it to its ports with strings.
Regions labelled by non-generating morphisms, such as identities, will have special rendering code.

This means that brick diagrams and string diagrams can be represented by the same underlying data structure.
It is even possible that the user could choose their preferred way to render the diagram.

\bibliographystyle{alpha}
\bibliography{\string~/Dropbox/Work/refs}

\end{document}